\documentclass[12pt]{article}

\usepackage{latexsym}
\usepackage{amssymb}
\usepackage{graphicx}
\usepackage{color}
\usepackage{amsmath}
\usepackage{float}
\usepackage{mathrsfs} 

\newtheorem{theorem}{Theorem}

\newtheorem{corollary}[theorem]{Corollary}

\newtheorem{definition}[theorem]{Definition}

\newtheorem{lemma}[theorem]{Lemma}

\newtheorem{remark}[theorem]{Remark}

\floatstyle{ruled} \restylefloat{figure} \restylefloat{table}

\newcommand{\beq}{\begin{equation}}
\newcommand{\bea}[1]{\begin{array}{#1} }
\newcommand{\eeq}{ \end{equation}}
\newcommand{\ea}{ \end{array}}

\newcommand{\ds}{\displaystyle}

\newcommand{\ar}{\partial}

\newcommand{\kom}[1]{}
\renewcommand{\kom}[1]{{\bf [#1]}}

\definecolor{blau}{rgb}{0.1,0.0,0.9}

\newcounter{komcounter}
\numberwithin{komcounter}{section}

\begin{document}

\title{
Phragm\'en-Lindel\"of theorems and $p$-harmonic measures for sets near low-dimensional hyperplanes}
\author{
Niklas L.P. Lundstr\"{o}m \linebreak \\\\
\it \small Department of Mathematics and Mathematical Statistics, Ume{\aa} University\\
\it \small SE-90187 Ume{\aa}, Sweden\/{\rm ;}
\it \small niklas.lundstrom@math.umu.se
}

\maketitle

\begin{abstract}
We prove estimates of a $p$-harmonic measure, $p \in (n-m, \infty]$,
for sets in $\mathbf{R}^n$ which are close to an $m$-dimensional hyperplane
$\Lambda \subset \mathbf{R}^n$, $m \in [0,n-1]$.
Using these estimates, we derive results of Phragm\'en-Lindel\"of type in unbounded domains $\Omega \subset \mathbf{R}^n\setminus \Lambda$ for $p$-subharmonic functions.
Moreover, we give local and global growth estimates for
$p$-harmonic functions,
vanishing on sets in $\mathbf{R}^n$,
which are close to an $m$-dimensional hyperplane. \\


\end{abstract}



\setcounter{equation}{0} \setcounter{theorem}{0}

\section{Introduction}


The $p$-harmonic functions,
which are natural nonlinear generalizations of the harmonic functions,
are solutions to the $p$-Laplace equation
\begin{align}\label{eq:plapequation}
\Delta_{p} u :=\nabla \cdot ( |\nabla u |^{p - 2} \nabla  u ) = 0,
\end{align}
when $p \in(1, \infty)$.
If $p = \infty$, then the equation can be written as
\begin{align}\label{eq:inflapequation}
\Delta_{\infty} u := \sum_{i,j = 1}^{n} \frac{\partial u}{\partial x_i} \frac{\partial u}{\partial x_j} \frac{\partial^2 u}{\partial x_i \partial x_j} = 0,
\end{align}
which is the so called $\infty$-Laplace equation.
We refer the reader to Section \ref{sec:prel} for the definitions of weak solutions,
viscosity solutions and $p$-harmonicity.
The $p$-Laplace equation has connections to e.g. minimization problems,
nonlinear elasticity theory,
Hele-Shaw flows and
image processing, see e.g. Lundstr\"om \cite[Chapter 2]{avhandlingen} and the references therein.


A class of $p$-harmonic functions that has shown to be useful consists of the following $p$-harmonic measures, 
which will be estimated in this paper.
\begin{definition} \label{def:p-hmeas}
Let $G \subseteq \mathbf{R}^n$ be a domain, $E \subseteq \partial{G}$, $p \in (1,\infty)$
and $x \in G$.
The $p$-harmonic measure of $E$ at $x$ with respect to $G$
is defined as $\inf_{u} u(x)$,
where the infimum is taken over all $p$-superharmonic functions $u \ge 0$ in $G$
such that $\liminf_{z \to y} u(z) \ge 1$, for all $y \in E$.
\end{definition}
\noindent
The $\infty$-harmonic measure is defined in a similar manner, but with $p$-superharmonicity replaced by absolutely minimizing, see Peres--Schramm--Sheffield--Wilson \cite[pages 173--174]{inf-tow}.
It turns out that the $p$-harmonic measure in Definition \ref{def:p-hmeas} is a
$p$-harmonic function in $\Omega$, bounded below by $0$ and bounded above by $1$.
For these and other basic properties of $p$-harmonic measure we refer the reader to
Heinonen--Kilpel\"ainen--Martio \cite[Chapter 11]{HKM}.
To avoid confusion,
we mention that there are at least three different $p$-harmonic measures in the literature.
Besides the $p$-harmonic measure above, we refer to the definitions given by
Bennewitz--Lewis \cite{BL2005} and Herron--Koskela \cite{HK97}.

The $p$-harmonic measure is useful when estimating solutions to the $p$-Laplace equation, 
see e.g. \cite[Theorem~11.9]{HKM}.
Recently, Lundstr\"om--Vasilis \cite{LV13} proved estimates for $p$-harmonic measures in the plane,
which, together with a result by Hirata \cite{H08}, yield properties of the $p$-Green function.
The $p$-harmonic measure is also useful when studying quasiregular mappings, see \cite[Chapter 14]{HKM}.
Moreover, the $p$-harmonic measure has a probabilistic interpretation
in terms of the zero-sum two-player game \emph{tug-of-war},
see Peres--Sheffield \cite{noisy-tow} and \cite{inf-tow},
in which also estimates for $p$-harmonic measure are proved, e.g. for porous sets.

Let $\Lambda \subset \mathbf{R}^n$ be an $m$-dimensional hyperplane, $m \in [0,n-1]$,
and introduce the notation
\begin{align}\label{eq:cylinder-def}
\Lambda_{s} = \left\{x \in \mathbf{R}^n : d(x,\Lambda) \leq s \right\}. 
\end{align}
Assume that $\Omega\subset \mathbf{R}^n$ is an unbounded domain 
with boundary $\partial \Omega$ close to $\Lambda$
in the sense that $\Lambda \subseteq \complement\Omega \subseteq \Lambda_s$ for some $s > 0$. 
Denote by $B(w,R)$ the open ball in $\mathbf{R}^n$ with center $w$ and radius $R$.
Suppose that $w \in \Lambda$, $p \in (n-m,\infty]$ and let $v_r$ be the $p$-harmonic measure of  
$\partial B(w,R) \setminus \complement \Omega$ at $x$ with respect to $B(w, R)\cap \Omega$.
In Theorem \ref{th:p-harmonic-measure} we prove that there exists a constant $C$ such that
\begin{align}\label{eq:simpli-res-phmeas}
\frac{1}{C} \leq v_R(x)\, R^\beta \leq C 
\end{align}
whenever $R$ is large enough and $\beta = (p-n+m)/(p-1)$ with $\beta = 1$ for $p = \infty$.

Next, we use this estimate to prove Corollary \ref{th:Phragmen-Lindelof},
which is an extended version of the classical result of Phragm\'en--Lindel\"of \cite{PL08}.
In particular, suppose that $u$ is $p$-subharmonic in an unbounded domain
$\Omega$ satisfying $\Omega \cap \Lambda = \emptyset$
and suppose that $\limsup_{z\to\partial\Omega} u(z) \leq 0$.
Then either $u \leq 0$ in the whole of $\Omega$ or it holds that
\begin{align}\label{eq:phragmen-1*}
\liminf_{R \to \infty} \left( \frac{1}{R^\beta} \ds\sup_{\partial B(w,R)\cap \Omega} u \right) > 0,
\end{align}
%
%
%
where $\beta$ is as in \eqref{eq:simpli-res-phmeas}.
When $\Omega = \mathbf{R}^n\setminus \Lambda_s$, the above growth rate is sharp.
Corollary \ref{th:Phragmen-Lindelof} generalizes a result of Lindqvist \cite{L85},
who studied the borderline case $p = n$, 
to hold in the exponent range $p \in (n-m,\infty]$.


The Phragm\'en-Lindel\"of principle, which has connections to elasticity theory, see e.g. Horgan \cite{Horgan}, 
Quintanilla \cite{Q93}, has been frequently studied during the last century.
To mention few papers, Ahlfors \cite{A37} extended results from \cite{PL08} to the upper half space of $\mathbf{R}^n$,
Gilbarg \cite{G52} and Serrin \cite{S54} considered more general elliptic equations of second order
and Vitolo \cite{V04} considered the problem in angular sectors.
Kurta \cite{K93} and Jin--Lancaster \cite{JL99, JL00, JL03} considered quasilinear elliptic equations and non-hyperbolic equations while
Capuzzo--Vitolo \cite{CV07} and Armstrong--Sirakov--Smart \cite{ASS12} considered fully nonlinear equations.
Adamowicz \cite{A14} studied different unbounded domains for subsolutions of the variable exponent $p$-Laplace equation,
while Bhattacharya \cite{Bhatt05} and Granlund--Marola \cite{GM14} considered infinity-harmonic functions in unbounded domains.

In connection with the above Phragm\'en-Lindel\"of result,
we also prove global growth estimates for positive $p$-harmonic functions,
vanishing on $\partial\Omega$,
where $\Omega$ is an unbounded domain as described above \eqref{eq:simpli-res-phmeas}.
This result is given in Theorem \ref{th:Phragmen-Lindelof-2} and implies,
in analogue with \eqref{eq:phragmen-1*}, that $u(x) \approx d(x,\partial\Omega)^\beta$ whenever $x \in \mathbf{R}^n$ and $d(x,\partial\Omega)$ is large.
Theorem \ref{th:Phragmen-Lindelof-2} generalizes e.g., some results by Kilpel\"ainen-Shahgholian-Zhong \cite{KSZ07} to hold in a more general geometric setting.

Our proofs rely on comparison with certain explicit $p$-subharmonic and $p$-superharmonic functions,
first constructed and used in Lundstr\"om \cite{Lu11}
to prove local estimates for $p$-harmonic functions. 
In this paper, we first expand this construction (Lemma \ref{le:superlosningen}),
through which we obtain an extension of all the main results in \cite{Lu11}, 
given for $p \in (n, \infty]$, 
to hold also in the wider exponent range $p \in (n-m,\infty]$ (Corollary \ref{cor:Lu11}).
Next, we use the explicit $p$-subharmonic and $p$-superharmonic functions in Lemma \ref{le:superlosningen} 
to prove local growth estimates (Theorem \ref{le:pipelemma}) 
for positive $p$-harmonic functions vanishing on a fraction of $\Lambda_s$.
The estimates in Theorem \ref{le:pipelemma}
are crucial for the proofs of our main results in Theorem \ref{th:p-harmonic-measure}, 
Corollary \ref{th:Phragmen-Lindelof} and Theorem \ref{th:Phragmen-Lindelof-2}.
Moreover, Theorem \ref{le:pipelemma} implies boundary Harnack's inequality near $\Lambda_s$ (Corollary \ref{cor:ny2:BHI}). 

Local estimates such as the boundary Harnack inequality for positive $p$-harmonic functions
vanishing on a portion of an $(n-1)$-dimensional boundary 
have drawn a lot of attention the last decades.
In the case $1 < p < \infty$, see e.g. 
Aikawa--Kilpel\"ainen--Shanmugalingam--Zhong~\cite{aikawa} for smooth boundaries,
Lewis--Nystr\"om~\cite{LN07, LN12} for more general geometries including Lipschitz and Reifenberg flat boundaries. 
For infinity-harmonic functions, see e.g. Bhattacharya~\cite{Bhatt}, Lundstr\"om--Nystr\"om~\cite{LuN} and
for solutions to the variable exponent $p$-Laplace equation in smooth domains, see Adamowicz--Lundstr\"om \cite{AL14}.
Only few papers considered local estimates of positive $p$-harmonic functions vanishing near boundaries having dimension less than $n-1$.
Besides results given in Theorem \ref{le:pipelemma} and Corollaries \ref{cor:ny1}--\ref{cor:ny2:BHI},
we refer the reader to Lindqvist \cite{L85} and Lundstr\"om \cite{Lu11}.
In an upcoming paper, Lewis and Nystr\"om will prove results in this direction for solutions to $p$-Laplace type equations near low-dimensional Reifenberg flat sets.


\section{Notation and preliminary lemmas}
\label{sec:prel}

By $\Omega$ we denote a domain, that is, an open connected set.
For a set $E \subset \mathbf{R}^n$
we let $\overline{E}$ denote the closure,
$\partial E$ the boundary and 
$\complement E$ the complement of $E$ and we put $E^o = E\setminus \partial E$.
Further, $d(x,E)$ denotes the Euclidean distance from $x \in \mathbf{R}^n$ to $E$, and
$B(x,r) = \{ y : | x  -  y | < r \}$
denotes the open ball with radius $r$ and center $x$.
By $c$ we denote a constant $\geq 1$, not necessarily the same at each occurrence,
depending only on $n$ and $p$ if nothing else is mentioned.
Moreover, $c(a_1, a_2, \dots, a_k)$ denotes a constant $\geq 1$, not necessarily the same at each occurrence,
depending only on $a_1, a_2, \dots, a_k$, and
we write $A \approx B$ if there exists a constant $c$ such that $c^{-1}A \leq B\leq cA$.
We denote points in Euclidean $n$-space $\mathbf R^{n}$ by
$x = (x_1,x_2,\dots,x_n) = (x',x'')$, where
\begin{align}\label{eq:notation-x'-and-x''}
x' = (x_1, x_2, \dots, x_{n-m}) \quad \textrm{and} \quad x'' = (x_{n-m+1}, x_{n-m+2},\dots, x_n).
\end{align}
Finally, we write $N = \{1,2,3,\dots\}$ for the set of natural numbers.

We next recall standard definitions of weak solutions, viscosity solutions and $p$-harmonicity.
If $p \in (1, \infty)$, we say that $u$ is a \textit{weak subsolution (supersolution)}
to the $p$-Laplace equation in $\Omega$ provided $u \in W_{loc}^{1,p}(\Omega)$ and
\begin{equation}\label{eq:huvudekva}
\int\limits_{\Omega} | \nabla  u |^{p - 2}  \, \langle   \nabla  u , \nabla \theta  \rangle \, dx \leq (\geq) \,0
\end{equation}
whenever $\theta \in C^{\infty}_0(\Omega)$ is non-negative. 
A function $u$ is a \textit{weak solution} of the
$p$-Laplacian if it is both a weak subsolution and a weak supersolution.
Here, as in the sequel,
$W^{1,p}(\Omega)$ is the Sobolev space of those $p$-integrable functions whose first
distributional derivatives are also $p$-integrable,
and $C^\infty_0(\Omega)$ is the set of infinitely differentiable functions with compact support in $\Omega$.
If $p = \infty$, the equation is no longer of divergence form and therefore the above definition is replaced by the definition of viscosity solutions, Crandall--Ishii--Lions \cite{guide}.

An upper semicontinuous function $u : \Omega \rightarrow \mathbf{R}$ is a \textit{viscosity subsolution} of the
$\infty$-Laplacian in $\Omega$ provided that for each function $\psi \in C^2(\Omega)$ such that $u-\psi$ has a local maximum at a point $x_0 \in \Omega$, we have $\Delta_{\infty}\psi(x_0) \geq 0$.
A lower semicontinuous function $u : \Omega \rightarrow \mathbf{R}$ is a \textit{viscosity supersolution} of the
$\infty$-Laplacian in $\Omega$ provided that for each function $\psi \in C^2(\Omega)$ such that $u-\psi$ has a local minimum at a point $x_0 \in \Omega$, we have $\Delta_{\infty}\psi(x_0) \leq 0$.
A function $u : \Omega \rightarrow \mathbf{R}$ is a \textit{viscosity solution} of the
$\infty$-Laplacian if it is both a viscosity subsolution and a viscosity supersolution.
In the following, we sometimes just write \textit{solution} for weak solutions and viscosity solutions.

If $u$ is an upper semicontinuous subsolution to the $p$-Laplacian in $\Omega$, $p \in (1, \infty]$, 
then we say that $u$ is \textit{$p$-subharmonic} in $\Omega$.
If $u$ is a lower semicontinuous supersolution to the $p$-Laplacian in $\Omega$, $p \in (1, \infty]$, 
then we say that $u$ is \textit{$p$-superharmonic} in $\Omega$.
If $u$ is a continuous solution to the $p$-Laplacian in $\Omega$, $p \in (1, \infty]$, then $u$ is
\textit{$p$-harmonic} in $\Omega$.

We note that for the $p$-Laplacian, $1 < p < \infty$, $p$-harmonic functions are equivalent to viscosity solutions 
(defined as above but with $\Delta_{\infty}$ replaced by $\Delta_p$); see Juutinen--Lindqvist--Manfredi \cite{JULM01}.
Moreover, in many situations, an $\infty$-harmonic function is the limit of a sequence of $p$-harmonic functions as $p \to \infty$; see Jensen \cite{J}.
This fact has been used to prove results for $p = \infty$ by taking limits of problems for finite $p$,
in which estimates are independent of $p$ when $p$ is large,
see e.g.
Bhattacharya--DiBenedetto--Manfredi \cite{BDM89},
Lindqvist--Manfredi \cite{lin-man},
Lewis--Nystr\"om~\cite{LN-bhi-plane} and Lundstr\"om--Nystr\"om~\cite{LuN}.
As for Phragm\'en-Lindel\"of type results,
see Granlund--Marola \cite{GM14}.
With this in mind, we chose to keep track of the dependence of $p$ in our estimates and point out
when constants are independent of $p$ when $p$ is large.


We next recall some well known results for $p$-harmonic functions.

\begin{lemma}\label{jamforelseprin}
(Comparison principle)
Let $p \in (1,\infty]$ be given, $u$ be $p$-subharmonic and $v$ be $p$-superharmonic in a bounded domain $\Omega$.
If
\begin{equation*}
\limsup_{x\to y} u(x) \leq \liminf_{x\to y} v(x)
\end{equation*}
for all $y \in \partial \Omega$, and if both sides of the above inequality are not simultaneously $\infty$ or $-\infty$,
then $u \leq v$ in $\Omega$.
\end{lemma}
\noindent
{\bf Proof.}
If $p \in (1, \infty)$, this result follows from Heinonen--Kilpel\"ainen--Martio \cite[Theorem~7.6]{HKM}.
For the case $p = \infty$, the lemma was first proved by Jensen \cite[Theorem 3.11]{J}.
Alternative proofs were later presented by Barles--Busca \cite{additional_ref} and Armstrong--Smart \cite{armstrongsmart}. $\hfill \Box$


\begin{lemma}\label{le:harnack}
(Harnack's inequality)
Let $p \in (1,\infty]$ be given and assume that $w \in \mathbf{R}^n$, $r \in (0, \infty)$ and that $u$ is a positive
$p$-harmonic function in $B(w,2r)$.
Then there exists $c(n,p)$, independent of $p$ if $p$ is large, such that
\begin{align*}
\sup_{B(w, r)} u \leq c \inf_{B(w, r)} u.
\end{align*}
\end{lemma}
\noindent
{\bf Proof.} For the case $p \in (1,\infty)$, when the constant is allowed to depend on $p$, we refer the reader to Heinonen--Kilpel\"ainen--Martio \cite[Theorem 6.2]{HKM}.
For the uniform in $p$ case, see Koskela--Manfredi--Villamor \cite{KMV}, Lindqvist--Manfredi \cite{lin-man} or
Lundstr\"om--Nystr\"om  \cite[Lemma~2.3]{LuN}.
For the case $p = \infty$ the result follows by taking the limit $p\to \infty$ in the above uniform in $p$ estimate;
see \cite{lin-man}. Moreover, another proof concerning the case $p = \infty$ was given by Bhattacharya \cite{Bhatt2}. $\hfill \Box$


\section{Estimates for $p$-harmonic functions vanishing near $m$-dimensional hyperplanes}
\label{sec:local}
\setcounter{equation}{0} \setcounter{theorem}{0}

We begin this section by stating, in our geometric setting,
some well known basic boundary estimates,
such as H\"older continuity up to the boundary and the Carleson estimate.
Next, we prove a refined version of Lundstr\"om \cite[Lemma 3.7]{Lu11} which yields explicit $p$-subharmonic and $p$-superharmonic functions, crucial for our proofs.
Moreover, we state and prove Theorem \ref{le:pipelemma}, giving growth estimates for $p$-harmonic functions 
vanishing near $m$-dimensional hyperplanes.
Finally, we discuss applications of Theorem \ref{le:pipelemma} by deriving several corollaries of the result.  

In the following we let $C_{p}$ denote $p$-capacity as defined in Heinonen--Kilpel\"ainen--Martio \cite[Chapter 2]{HKM}.
That is, 
the $p$-capacity of the condenser $(K, \Omega)$,
where $K\subset\Omega$ is compact,
is the number defined by
%
$$
C_p(K, \Omega) = \inf_u \int_{\Omega} |\nabla u|^p dx,
$$
where the infimum is taken over all $u \in C_0^{\infty}(\Omega)$ such that $u \geq 1$ on $K$.
%
%
%
%

\begin{lemma}\label{le:capacity-Manifold}
Let $M \subset \mathbf{R}^n$ be a manifold of dimension $m < n$,
then $M$ has $p$-capacity zero if and only if $p \leq n - m$.
\end{lemma}

\noindent
{\bf Proof.} The result follows from Adams--Hedberg \cite[Corollary 5.1.15]{AH96}. $\hfill\Box$


\begin{lemma}\label{le:holder}
(H\"older continuity)
Suppose that $m, n \in N$ such that $m \in [0,n-1]$,
let $\Lambda \subset \mathbf{R}^n$ be an $m$-dimensional hyperplane, $w \in \Lambda$, $r \in (0,\infty)$ and $p \in (n-m,\infty]$.
Assume that $u$ is a non-negative $p$-harmonic function in $B(w, 2r) \setminus \Lambda$,
continuous in $B(w, 2r)$ with $u = 0$ on $B(w, 2r) \cap \Lambda$.
Then there exist constants $\gamma \in (0, 1]$ and
$c$, both depending only on $p$ and $n$, independent of $p$ if $p$ is large,
such that if $x, y \in B(w, r)$ then
\begin{align*}
|u(x) - u(y)| \leq c \left( \frac{|x - y|}{r} \right)^{\gamma} \sup_{B(w, 2r)} u.
\end{align*}
In particular, we can take $\gamma \to 1$ as $p \to \infty$ with $\gamma = 1$ if $p = \infty$.
\end{lemma}

\noindent
{\bf Proof.}
If $p > n$ we obtain the result by a Sobolev embedding theorem, see e.g., Lundstr\"om--Nystr\"om  \cite[Lemma 2.4]{LuN}.
If $n-m < p \leq n$, then the lemma follows from Heinonen--Kilpel\"ainen--Martio \cite[Theorem 6.44]{HKM}
if we can prove that there exist constants $c_0$ and $r_0$ so that
\begin{align}\label{eq:cap-estimate}
\frac{C_p\left( \Lambda \cap \overline{B}(x_0,r), B(x_0, 2 r) \right)}{C_p \left( \overline{B}(x_0,r), B(x_0, 2r) \right)} \geq c_0
\end{align}
whenever $0 < r < r_0$ and $x_0 \in \Lambda$.
To prove \eqref{eq:cap-estimate} observe that,
since the $p$-capacity is invariant through rotations and translations,
it holds that
\begin{align*}
C_p\left( \Lambda \cap \overline{B}(x_0,r), B(x_0, 2 r) \right) =
r^{n-p} C_p\left( \{x\in\mathbf{R}^n : |x'| = 0\} \cap \overline{B}(0,1), B(0, 2) \right).
\end{align*}
Moreover, from Lemma \ref{le:capacity-Manifold} it follows, since $n-m < p$, that there exits $c(n,p)$ such that
\begin{align*}
C_p\left(  \{x\in\mathbf{R}^n : |x'| = 0\} \cap \overline{B}(0,1), B(0, 2) \right) \geq c^{-1} > 0.
\end{align*}
Since \cite[Example 2.12]{HKM} gives
$C_p \left( \overline{B}(x_0,r), B(x_0, 2r) \right) = c(n,p) r^{n-p}$,
inequality \eqref{eq:cap-estimate} follows for $r_0 = \infty$.
The proof of Lemma \ref{le:holder} is complete. $\hfill\Box$\\


Given an $m$-dimensional hyperplane $\Lambda$ and $w \in \Lambda$
we let in the following $A_r(w)$ denote a point satisfying
\begin{align}\label{eq:point_A_r(w)}
 d(A_r(w),\Lambda) = r \quad \textrm{and} \quad A_r(w) \in \partial B(w,r).
\end{align}
\begin{lemma}\label{le:carleson}
(Carleson's estimate)
Suppose that $m, n \in N$ such that $m \in [0,n-1]$,
let $\Lambda \subset \mathbf{R}^n$  be an $m$-dimensional hyperplane, $w \in \Lambda$, $r \in (0,\infty)$ and $p \in (n-m,\infty]$.
Assume that $u$ is a non-negative $p$-harmonic function in $B(w, r) \setminus \Lambda$,
continuous in $B(w, r)$ with $u = 0$ on $B(w, r) \cap \Lambda$.
Then there exists $c(n,p)$, independent of $p$ if $p$ is large, such that
\begin{eqnarray*}
\sup_{B(w, r/c)} u  \, \leq  \, c \,u(A_{r/c}(w)).
\end{eqnarray*}
\end{lemma}

\noindent
{\bf Proof.}
A proof for linear elliptic partial differential equations,
in Lipschitz domains with $(n-1)$-dimensional boundary,
can be found in Caffarelli--Fabes--Mortola--Salsa \cite{CFMS}.
The proof uses only the Harnack chain condition (see e.g. \cite[Definition 1.3]{lower_order}),
analogues of Harnack's inequality,
H\"older continuity up to the boundary
and the comparison principle for linear equations.
In particular, the proof also applies in our situation. $\hfill\Box$\\


The following lemma extends constructions in Lundstr\"om \cite[Lemma 3.7]{Lu11}, given for $p \in (n, \infty)$,
to hold for the wider exponent range $p \in (n-m, \infty)$.
Recall from \eqref{eq:notation-x'-and-x''} the notation $x = (x',x'')\in \mathbf{R}^n$ and
the geometric definition of $\Lambda_s$ given in \eqref{eq:cylinder-def} as
\begin{align*}
\Lambda_{s} = \left\{x \in \mathbf{R}^n : d(x,\Lambda) \leq s \right\}
\end{align*}
where $\Lambda$ is an $m$-dimensional hyperplane.

\begin{lemma}\label{le:superlosningen}
Suppose that $m, n \in N$ such that $m \in [1,n-2]$.
Let $p \in (n-m,\infty)$, $\beta = (p-n+m)/(p-1)$ and suppose that $\gamma$ satisfies $0 < \gamma < \beta$.
Then there exists $\delta_c \in (0,1)$, depending only on $n, \gamma$ and $p$,
such that $\hat u$ is a supersolution, and $\check u$ is a subsolution to the $p$-Laplace equation in
$\left\{ x : |x''| < 1\right\} \cap \Lambda^o_{\delta_c}\setminus \Lambda$, where
\begin{align*}
\hat u = |x'|^{\beta} + |x''|^2 |x'|^{\gamma} - \frac12 |x'|^2 \quad \textrm{and} \quad
\check u = (1 - |x''|^2) |x'|^{\beta} + |x'|.
\end{align*}
Moreover, if $\gamma > 1/2$ then $\delta_c$ can be chosen independent of $p$ if $p$ is large.
\end{lemma}
\noindent
{\bf Proof.}
For a proof showing that $\check u$ is a subsolution,
as well as for the case $\gamma = (p-n)/(p-1)$,
we refer the reader to the proof of Lemma 3.7 in \cite{Lu11}.
It remains to show that $\hat u$ is a supersolution for any $\gamma$,
$0 < \gamma < \beta$.
To do so, it suffices to show that there exists $\delta_c \in (0,1)$,
depending only on $\gamma, n$ and $p$,
such that
\begin{align}\label{eq:plaplaceutskriven1}
\Delta_p\hat u = \Delta \hat u |\nabla \hat u|^{p - 2}  +  (p - 2)|\nabla \hat u|^{p - 4}\Delta_{\infty}\hat u  \leq 0  \quad \textrm{in} \quad \left\{ x : |x''| < 1\right\} \cap \Lambda^o_{\delta_c}\setminus \Lambda.
\end{align}

\noindent
Here, $\Delta_p$ is the $p$-Laplace operator defined in \eqref{eq:plapequation}, $\Delta := \Delta_2$ and $\Delta_{\infty}$ is the $\infty$-Laplace operator defined in \eqref{eq:inflapequation}.
Since $p > n - m \geq 2$ and $|\nabla \hat u| \neq 0$ outside of $\Lambda$,
\eqref{eq:plaplaceutskriven1} equals
\begin{align}\label{eq:plaplaceutskriven}
\widehat\Delta_p\hat u := \frac{\Delta \hat u |\nabla \hat u|^2}{p - 2} + \Delta_{\infty}\hat u \leq 0 \quad \textrm{in} \quad \left\{ x : |x''| < 1\right\} \cap \Lambda^o_{\delta_c}\setminus \Lambda.
\end{align}

\noindent
Following the calculations in \cite[Pages 6857--6858]{Lu11} we obtain that
\begin{align}\label{eq:extraexplanation}
\widehat\Delta_p\hat u = \frac{\Delta \hat u |\nabla \hat u|^2}{p - 2} + \Delta_{\infty}\hat u
= Z_0 + Z_2 |x''|^2 + Z_4 |x''|^4 + Z_6 |x''|^6,
\end{align}

\noindent
where the coefficients are given by
\begin{align}\label{eq:ZZZZ}
Z_0 &= -\frac{\beta^2(p + n - 2 - m)}{p-2}|x'|^{2\beta - 2} +  O(|x'|^{\gamma + 2\beta - 2}) \nonumber \\
&\leq - \beta^2|x'|^{2\beta - 2}  +  O(|x'|^{\gamma + 2\beta - 2}),\nonumber\\
Z_2 &= - Z  \frac{\gamma\beta^2}{p-2}|x'|^{\gamma + 2\beta - 4}  
 -z_2 \gamma \beta |x'|^{\gamma + \beta -2}
 + O(|x'|^{2\gamma + \beta - 2}),\nonumber\\
Z_4 &= - Z \frac{2\gamma^2\beta}{p-2}|x'|^{2\gamma + \beta - 4}  
 -z_4 \gamma^2 |x'|^{2\gamma - 2}  +  O(|x'|^{3\gamma - 2}),\\
Z_6 &= - Z \frac{\gamma^3}{p-2}|x'|^{3\gamma - 4},\nonumber
\end{align}

\noindent
in which $Z = p - n + m - (p-1) \gamma$.
Clearly $Z > 0$ by the assumption $0 < \gamma < \beta$ and, hence,
%
%
we conclude that the leading terms are negative in \eqref{eq:ZZZZ}.
It follows from \eqref{eq:extraexplanation} and \eqref{eq:ZZZZ} that there exists $\delta_c \in (0,1)$,
depending only on $\gamma, n$ and $p$,
such that \eqref{eq:plaplaceutskriven} is satisfied in
$\left\{ x : |x''| < 1\right\} \cap \Lambda^o_{\delta_c}\setminus \Lambda$.

For the uniform in $p$ case, we note that if $p$ is large enough, then
\begin{align}\label{eq:johej-indep-of-p}
z_2 &= \frac{2 \gamma (p^2 - 2p + 1) + (3p - 2)(n - m - 1)}{p^2 - 3p + 2} \geq 2 \gamma, \notag\\
z_4 &= \frac{2 \gamma (p - 1) - p - 2 + 3(n - m)}{p-2} \geq 2 \gamma - 1.
\end{align}
By following calculations in \cite[Pages 6857--6858]{Lu11},
we see that the constants in the Ordos in \eqref{eq:ZZZZ} will not explode as $p \to \infty$.
Therefore, from \eqref{eq:ZZZZ}, \eqref{eq:johej-indep-of-p} and the assumption $\gamma > 1/2$,
we conclude that $\delta_c$ can be chosen independent of $p$ if $p$ is large,
but still depending on $n$ and $\gamma$.
This completes the proof of Lemma \ref{le:superlosningen}.
$\hfill \Box$\\


We are now ready to state and prove the main theorem of this section,
which gives the following upper and lower growth estimates of $p$-harmonic functions,
$p \in (n-m, \infty]$, vanishing near an $m$-dimensional hyperplane $\Lambda$. 
Recall the definition of $A_r(w)$ given in \eqref{eq:point_A_r(w)}. 

\begin{theorem}\label{le:pipelemma}
Suppose that $m, n \in N$ such that $m \in [0,n-1]$,
let $\Lambda \subset \mathbf{R}^n$ be an $m$-dimensional hyperplane,
$w \in \Lambda$, $r \in (0,\infty)$, $p \in (n-m,\infty]$
and suppose that $\beta = (p-n+m)/(p-1)$ with $\beta = 1$ if $p = \infty$.
Let $\delta \in [0, \delta_c / 2)$ where $\delta_c$ is from Lemma \ref{le:superlosningen} and
assume that $u$ is a positive $p$-harmonic function in $B(w,4r) \setminus \Lambda_{\delta r}$,
with $u = 0$ continuously on $B(w,4r) \cap \partial \Lambda_{\delta r}$.
Then there exists $c(n,p)$, independent of $p$ if $p$ is large, such that
\begin{align*}
c^{-1} \, \left\{\left(\frac{d(x, \Lambda)}{r}\right)^{\beta} - \delta^{\beta} \right\}
\leq \frac{u(x)}{u(A_{r}(w))}
\leq c \, \left\{\left(\frac{d(x, \Lambda)}{r}\right)^{\beta} - \delta^{\beta} \right\}
\end{align*}
whenever $x \in B(w,\delta_c\, r) \setminus \Lambda_{\delta r}$.
\end{theorem}

\noindent
Before proving the theorem, we make some remarks about the result.
For any $\delta \in (0, \delta_c / 2)$,
Theorem \ref{le:pipelemma} implies that, close to $\Lambda_{\delta r}$,
the $p$-harmonic function $u$ vanishes at the same rate as the distance function,
$u(x) \approx d(x,\Lambda)$, with constants exploding as $\delta \to 0$ unless $p = \infty$.
In particular, we obtain the following.

\begin{corollary}\label{cor:ny1}
Suppose that $m,n, \Lambda, w, r, p, \beta, \delta$ and $u$ are as in Theorem \ref{le:pipelemma}.
If $\delta > 0$,
then there exists $c(n,p)$, independent of $p$ if $p$ is large, such that
\begin{align*}
c^{-1} \, \delta^{\beta -1}  \frac{d(x, \Lambda_{\delta r})}{r}
\leq \frac{u(x)}{u(A_{r}(w))}
\leq c \, \delta^{\beta -1}  \frac{d(x, \Lambda_{\delta r})}{r}
\end{align*}
whenever $x \in B(w,2\delta r) \setminus \Lambda_{\delta r}$. 
\end{corollary}

\noindent
{\bf Proof.} 
The result follows by Taylor-expanding the estimates in Theorem \ref{le:pipelemma}. $\hfill \Box$\\

If $\delta = 0$ in Theorem \ref{le:pipelemma},
then $u \approx d(x, \Lambda)^{\beta}$.
In fact, for $\delta = 0$ we obtain the following corollary,
in which $C^{0,\beta}(E)$ denotes the space of H\"older continuous
functions in $E \subset \mathbf{R}^n$.

\begin{corollary}\label{cor:Lu11}
Suppose that $m,n, \Lambda, w, r, p, \beta, \delta$ and $u$ are as in Theorem \ref{le:pipelemma}.
If $\delta = 0$,
then 
%
\begin{align}\label{eq:first-in-corollary}
c^{-1} \, \left(\frac{d(x, \Lambda)}{r}\right)^{\beta}
\leq \frac{u(x)}{u(A_{r}(w))}
\leq c \, \left(\frac{d(x, \Lambda)}{r}\right)^{\beta}
\end{align}
whenever $x \in B(w,\delta_c\, r) \setminus \Lambda$ and $c$ is the constant from Theorem \ref{le:pipelemma}. 
Moreover, there exists $c(n,p)$ such that $u \in C^{0,\beta}(B(w,r/c))$,
and $\beta$ is the optimal H\"older exponent for $u$.
\end{corollary}

\noindent
{\bf Proof.} 
Estimate \eqref{eq:first-in-corollary} follows immediately by taking $\delta = 0$ in Theorem \ref{le:pipelemma}.
Using \eqref{eq:first-in-corollary} in place of \cite[Theorem 1.1]{Lu11},
and observing from Kilpel\"ainen--Zhong \cite{KZ, KZ03} that Lemma 2.4 in \cite{Lu11}
holds also in the wider exponent range $p \in (n-m, \infty)$,
the H\"older continuity follows by mimicking the proof of Corollary 1.2 in \cite{Lu11}. $\hfill \Box$\\

\noindent
Corollary \ref{cor:Lu11} retrieves the geometric setting of \cite[Theorem 1.1 and Corollary 1.2]{Lu11},
and generalizes these theorems, given for $p \in (n, \infty]$,
to hold also in the wider exponent range $p \in (n-m, \infty]$.

Moreover, since Theorem \ref{le:pipelemma} gives the rate at which $p$-harmonic functions vanish near $\Lambda_{\delta r}$,  it implies the boundary Harnack inequality: 

\begin{corollary}\label{cor:ny2:BHI}
(boundary Harnack's inequality) Suppose that $m,n, \Lambda, w, r, p, \beta, \delta$ and $u$ are as in Theorem \ref{le:pipelemma}.
Assume that $v$ is a $p$-harmonic function satisfying the same assumptions as $u$, then
$$
c^{-2}\, \frac{u(A_r(w))}{v(A_r(w))}
\,\leq
\,\frac{u(x)}{v(x)}
\,\leq
\,c^2\, \frac{u(A_r(w))}{v(A_r(w))} 
$$
whenever $x \in B(w,\delta_c r) \setminus \Lambda_{\delta r}$ and $c$ is the constant from Theorem \ref{le:pipelemma}. 
\end{corollary}

\noindent
{\bf Proof.} 
The corollary follows by applying Theorem \ref{le:pipelemma} to the $p$-harmonic functions $u$ and $v$. $\hfill \Box$\\

Besides the applications above and those given in Section \ref{sec:p-harmonic-measures},
Theorem \ref{le:pipelemma} can be useful when studying local estimates of
$p$-harmonic functions vanishing on sets which can be trapped into $\Lambda_s$.
An example of such sets are the $m$-dimensional Reifenberg-flat sets,
which are approximable, uniformly on small scales, by $m$-dimensional hyperplanes.
For the definition of Reifenberg-flat sets and for some applications, involving boundary behaviour of solutions to PDEs,
see e.g. 
Kenig--Toro \cite{KT97},
David \cite{D14},
Guanghao--Wang \cite{GW07},
Capogna--Kenig--Lanzani \cite{CKL05},
Lewis--Nystr\"om \cite{LN12} 
and Avelin--Lundstr\"om--Nystr\"om \cite{lower_order, optimal-doubling}. \\





\noindent
{\bf Proof of Theorem \ref{le:pipelemma}.}
Since the $p$-Laplace equation is invariant under scalings, translations and rotations,
we assume, without loss of generality,
that $w=0$, $r = 1$, $u(A_r(w)) = u(A_1(0)) = 1$ and
\begin{align*}
\Lambda = \{x \in \mathbf{R}^n : |x'| = 0\}.
\end{align*}
In these coordinates, we will prove
the existence of $c(n,p)$ such that
\begin{align}\label{eq:scaled_simple_0}
c^{-1} \, \left\{  |x'|^{\beta} - \delta^{\beta} \right\}
\leq u(x)
\leq c \, \left\{  |x'|^{\beta} - \delta^{\beta} \right\}
\end{align}
whenever $x \in B(0,\delta_c) \setminus \Lambda_{\delta}$.
Scaling back then yields Theorem \ref{le:pipelemma}. \\

\noindent
{\bf Proof of the upper bound.}
%
%
We begin with the case $m = n-1$, in which
the Theorem follows by already well known results,
such as e.g. Aikawa--Kilpel\"ainen--Shanmugalingam--Zhong \cite{aikawa}.
We include a proof for the sake of completeness.
Since, in this case, $\Lambda$ splits $\mathbf{R}^n$ in two halves,
we focus on the upper of these halves.
Let $\alpha = (p-n)/(p-1)$ with $\alpha = 1$ if $p = \infty$ and
consider the $p$-harmonic function
\begin{eqnarray}\label{eq:first-fundamental-f}
\bar f(x) &=& a |x - x_0|^{\alpha} + b, \quad\qquad \textrm{if} \qquad p \neq n,  \notag \\
\bar f(x) &=& a \log|x - x_0| + b, \qquad \textrm{if} \qquad p = n,
\end{eqnarray}
for some $a,b$.
Choose $a$ and $b$ such that $\bar f$ has boundary values $\bar f = 0$ on $\ar B(x_0, 1/2)$
and $\bar f = 1$ on $\ar B(x_0, 1)$.
From \eqref{eq:first-fundamental-f} we conclude the existence of $c(n,p)$, decreasing in $p$,
such that
\begin{align}\label{eq:fundam-grad-beg}
c^{-1} \leq \frac{ \partial \bar f}{\partial \nu} \leq c
\quad \textrm{in} \quad \overline{B(x_0, 1)} \setminus B(x_0,1/2),
\end{align}
where $\nu$ denotes the outer normal to $\partial B(x_0, 1)$.
Since $u(A_1(0)) = 1$ there exists, by Harnack's inequality and the Carleson estimate,
a constant $\bar c(n,p)$ such that
\begin{align*}
u(x) \leq \bar c \quad \textrm{in} \quad B(0,3) \cap \{x : x' = x_1 > \delta\}.
\end{align*}
Since $u$ vanishes continuously on $\partial \Lambda_{\delta} \cap B(0,4)$,
we can conclude,
by the comparison principle applied to the functions
$u$ and $\bar c \bar f$
and by letting $x_0$ vary with the restriction that $B(x_0,1/2)$ is tangent to $\{x : x_1 = \delta\}$,
$B(x_0,1/2) \subset \{x : x_1 < \delta\}$ and $B(x_0,1) \subset B(0,3)$,
that there exists $c(n,p)$, independent of $\delta$ and independent of $p$ if $p$ is large, such that
\begin{align*}
u(x) \leq c \left( |x'| - \delta \right) \quad \textrm{whenever}\quad x \in B(0,1) \cap \{x : x' = x_1 > \delta\}.
\end{align*}
Thus, we have proved the upper bound in Theorem \ref{le:pipelemma} in the case $m = n-1$.

%
In the rest of the proof of the upper bound, we assume $m \in [0,n-2]$.
%
%
Assume first also that $p > n$ and consider the $p$-harmonic function
\begin{align}\label{eq:upper-f-hat}
\hat f(x) =  |x - x_0|^{\alpha} - \delta^{\alpha},
\end{align}
where $x_0 \in {\Lambda}\cap B(0,2)$ and $\alpha$ is the exponent defined above \eqref{eq:first-fundamental-f}.
Note that $\hat f \geq 0$ on $B(x_0,1) \setminus \Lambda_\delta$ and $\hat f = 1-\delta^{\alpha}$ on
$\partial B(x_0,1)$.
Using $u(A_1(0)) = 1$, Harnack's inequality and the Carleson estimate,
we will now show that there exists a constant $\hat c(n,p)$,
independent of $\delta$, such that
\begin{align}\label{eq:ny-simpel-men-inte-trivial}
u(x) \leq \hat c \quad \textrm{in} \quad B(0,3)\setminus \Lambda_\delta.
\end{align}
To prove \eqref{eq:ny-simpel-men-inte-trivial},
let $\tilde u$ be the $p$-harmonic function in e.g. $B(0,3\frac{1}{2})\setminus \Lambda$,
satisfying boundary values
$\tilde u = u$ on $\partial B(0,3\frac{1}{2})\setminus \Lambda_\delta$ and
$\tilde{u} = 0$ on $\Lambda \cup (\partial B(0,3\frac{1}{2})\cap \Lambda_\delta)$ continuously.
Note that the boundary values for $\tilde u$ are continuous and that existence of
$\tilde u$ follows from \eqref{eq:cap-estimate} and standard existence theorems,
see e.g. Heinonen--Kilpel\"ainen--Martio \cite{HKM}.
It follows by construction and by the comparison principle that
$u \leq \tilde u$ in $B(0,3\frac{1}{2})\setminus \Lambda_\delta$. 
Applying Harnack's inequality and the Carleson estimate to $\tilde u$ implies \eqref{eq:ny-simpel-men-inte-trivial}.
Since $u$ vanishes continuously on $\partial \Lambda_{\delta} \cap B(0,2)$
and $\hat f = 0$ on $B(x_0, \delta)$,
it follows by the comparison principle, applied to $u$ and $\hat c \hat f$
and by letting $x_0 \in \Lambda \cap B(0,2)$ vary,
that there exists $c(n,p)$, independent of $\delta$ and independent of $p$ if $p$ is large, such that
\begin{align}\label{eq:upp1}
u(x) \leq c \left( |x'|^{\alpha} - \delta^{\alpha}\right) \quad
\textrm{whenever}\quad x \in B(0,2) \setminus \Lambda_\delta.
\end{align}
If $p = \infty$ or if $m = 0$, then
we have proved the upper bound in Theorem \ref{le:pipelemma}.
%
%

We now assume $n - m < p \leq n$ (implying $m \geq 1$) and prove that,
by H\"older continuity up to the boundary, $u(A_1(0)) = 1$, Harnack's inequality and the Carleson estimate,
there exist $c(n,p)$ and $\bar \gamma(n,p)$,
independent of $\delta$ and independent of $p$ if $p$ is large,
such that
\begin{align}\label{eq:upp1-med-n-m<p<=n}
u(x) \leq c \, |x'|^{\bar \gamma}  \quad \textrm{whenever}\quad x \in B(0,2) \setminus \Lambda_\delta.
\end{align}
To prove \eqref{eq:upp1-med-n-m<p<=n},
consider the auxiliary function $\tilde u$ defined below \eqref{eq:ny-simpel-men-inte-trivial}
but with $\tilde{u} = 0$ on $(\partial B(0,3\frac{1}{2})\cap \Lambda_\delta) \cup \widetilde\Lambda$,
instead of $\tilde u = 0$ on $\Lambda$, where $\widetilde\Lambda$ is an $m$-dimensional hyperplane
parallel to $\Lambda$ satisfying $\widetilde{\Lambda} \subset \Lambda_\delta$.
As before, it follows that $u \leq \tilde{u}$ in $B(0,3\frac{1}{2})\setminus \Lambda_\delta$.
Allowing $\widetilde\Lambda$ to move in $\Lambda_\delta$
and by using Lemma \ref{le:holder} (H\"older continuity), the Carleson estimate, Harnack's inequality and $\tilde u(A_1(0)) \approx u(A_1(0)) = 1$,
we conclude \eqref{eq:upp1-med-n-m<p<=n}.

Using estimates \eqref{eq:upp1} and \eqref{eq:upp1-med-n-m<p<=n} we will now use the supersolution given in
Lemma \ref{le:superlosningen} to complete the proof of the upper bound for the remaining cases
$m \in [1,n-2]$ and $p \in (n-m,\infty)$.
To do so, we will first show that there exists $c$ such that
%
\begin{align}\label{eq:domain-for-supersolution}
u \leq c \left(\hat u - \delta^{\beta} + \frac12\delta^{2}\right) \quad \textrm{on} \quad
\partial (\{x : |x''| \leq 1\} \cap \Lambda^o_{\delta_c} \setminus \Lambda_{\delta}).
\end{align}
Recall the assumption $2\delta < \delta_c < 1$. Using the definition of $\hat u$ it follows that on this set we have either
\begin{align}\label{eq:svans}
&|x'| = \delta, \quad \textrm{implying} \quad \hat u - \delta^{\beta} + \frac12\delta^{2}
\geq 0, \quad \mbox{or}\\
&|x'| = \delta_c, \quad \textrm{implying} \quad \hat u - \delta^{\beta} + \frac12\delta^{2}
\geq \delta_c^{\beta} - \delta^{\beta} - \frac12 \delta_c^2 + \frac12\delta^{2}
\geq \frac1{c}, \quad \mbox{or}\nonumber\\
&|x''| = 1\; \textrm{and}\; \delta < |x'| < \delta_c, \quad \textrm{implying}\nonumber\\
& \hat u - \delta^{\beta} + \frac12\delta^{2}
= |x'|^{\beta} + |x'|^{\gamma} - \frac12|x'|^2 - \delta^{\beta} + \frac12\delta^{2}
\geq |x'|^{\gamma} - \frac12 |x'|^2 + \frac12\delta^{2}
\geq \frac12 |x'|^{\gamma},\nonumber
\end{align}
for some $c$ depending only on $\beta$ and $\delta_c$.
From \eqref{eq:upp1}, \eqref{eq:upp1-med-n-m<p<=n} and \eqref{eq:svans} we conclude
\eqref{eq:domain-for-supersolution} by taking $\gamma = \alpha$ or $\gamma = \bar \gamma$ in \eqref{eq:svans}.
By the the comparison principle we obtain
\begin{align*}
u \leq c \left(\hat u - \delta^{\beta} + \frac12\delta^{2}\right) \quad \textrm{in} \quad
\{x : |x''| \leq 1\} \cap \Lambda^o_{\delta_c} \setminus \Lambda_{\delta}.
\end{align*}
By the definition of $\hat u$ it follows that
\begin{align}\label{eq:sistaupp}
u(x)
\leq c \left(|x'|^{\beta} + |x''|^2 |x'|^{\gamma} - \frac12 |x'|^2 - \delta^{\beta} + \frac12\delta^{2}\right)
\leq c \left(|x'|^{\beta} - \delta^{\beta}\right)
\end{align}
whenever 
$x \in \{x : |x''| = 0\} \cap \Lambda^o_{\delta_c} \setminus \Lambda_{\delta}$.
The constants in \eqref{eq:sistaupp} depend only on $n, p$ and $\delta_c$,
where $\delta_c(n,p,\gamma)$ is from Lemma \ref{le:superlosningen}.
Since, by Lemma \ref{le:holder},
$\gamma = \gamma(n,p)$ and $\gamma \to 1$ as $p \to \infty$,
we conclude, from Lemma \ref{le:superlosningen},
that the constants in \eqref{eq:sistaupp} depend only on $n, p$, independent of $p$ if $p$ is large.

Finally, by translating the function $\hat u - \delta^{\beta} + \frac{1}{2}\delta^2$
 and the domain 
$\{x : |x''| \leq 1\} \cap \Lambda^o_{\delta_c} \setminus \Lambda_{\delta}$
in the $x''$-direction, we finish the proof of the upper bound.
In particular, as long as
$\{x : |x''| \leq 1\} \cap \Lambda^o_{\delta_c} \setminus \Lambda_{\delta} \subset B(0,2)$
where we have \eqref{eq:upp1} and \eqref{eq:upp1-med-n-m<p<=n},
we may apply the same argument.
Thus we obtain that \eqref{eq:sistaupp} holds true in $B(0,\delta_c)$,
which completes the proof of the upper bound in \eqref{eq:scaled_simple_0} and hence also in Theorem \ref{le:pipelemma}.\\

\noindent
{\bf Proof of the lower bound.}
We first observe that since $u(A_1(0)) = 1$ we obtain by Harnack's inequality
(focusing on the upper half of $\mathbf{R}^n$ when $m = n-1$) that
\begin{align}\label{eq:johej_ny_revise_1}
 c^{-1} \leq u(x) \quad \textrm{on} \quad  \mbox{$B(0,3\frac{1}{2})$}\setminus \Lambda_{1/2 + \delta}.
\end{align}
If $m = 0$, then we use comparison with the function $\hat f$ from \eqref{eq:upper-f-hat} as follows.
Put $x_0 = 0$ and observe that then $c^{-1}\hat f \leq u$ on $\partial B(0,1)\cup\partial \Lambda_\delta$.
By the comparison principle $c^{-1}\hat f \leq u$ in $B(0,1)\setminus\Lambda_\delta$ and so
\begin{align*}
|x|^{\alpha} - \delta^\alpha \leq c \, u(x) \quad \textrm{whenever} \quad x \in B(0,1) \setminus \Lambda_{\delta}.
\end{align*}
This proves the lower bound when $m = 0$.

Next, assume that $m \geq 1$ and consider the $p$-harmonic function
\begin{eqnarray*}
\check f(x) &=& a |x - x_0|^{\alpha} + b, \quad\qquad \textrm{if} \qquad p \neq n,  \notag \\
\check f(x) &=& a \log|x - x_0| + b, \qquad \textrm{if} \qquad p = n,
\end{eqnarray*}
for some $a,b$ and with $\alpha$ defined as above \eqref{eq:first-fundamental-f}.
Choose $a$ and $b$ such that $\check f$ has boundary values
$\check f = 0$ at $\ar B(x_0, 1)$ and $\check f = 1$ at $\ar B(x_0, 1/2)$.
Using 
\eqref{eq:fundam-grad-beg} we see that
$c \check f \geq 1 - |x-x_0|$ in $B(x_0, 1)\setminus B(x_0, 1/2)$ for some $c(n,p)$ decreasing in $p$.
Using \eqref{eq:johej_ny_revise_1} and $c^{-1}\check f$ as a barrier from below for $u$ by placing the ball $B(x_0, 1)$ tangent to $\Lambda_{\delta}$
and allowing $x_0$ to vary, with the restriction $B(x_0, 1) \subset B(0,3\frac{1}{2})$,
we see that there exists $c(n,p)$ such that
\begin{align}\label{eq:ner1}
|x'| - \delta \leq c \, u(x) \quad \textrm{whenever} \quad x \in B(0,2) \setminus \Lambda_{\delta}.
\end{align}

\noindent
If $m = n - 1$ or if $p = \infty$, then the lower bound in Theorem \ref{le:pipelemma} follows from \eqref{eq:ner1}.

We assume from now on that $m \in [1, n-2]$ and $p \in (n-m,\infty)$.
The next step is to use the subsolution $\check u - (\delta^\beta + \delta)$,
derived in Lemma \ref{le:superlosningen}, as follows.
On $\partial (\{x : |x''| \leq 1\} \cap \Lambda^o_{\delta_c} \setminus \Lambda_{\delta})$ we have either
\begin{align}\label{eq:svans2}
&|x'| = \delta, \quad \textrm{implying} \quad \check u - (\delta^\beta + \delta)
= (1-|x''|^2)\delta^{\beta} + \delta - (\delta^\beta + \delta) \leq 0, \quad \mbox{or} \nonumber\\
&|x'| = \delta_c, \quad \textrm{implying} \quad \check u - (\delta^\beta + \delta)
= (1-|x''|^2)\delta_c^{\beta} + \delta_c - (\delta^\beta + \delta)
\leq c, \quad \mbox{or} \nonumber\\
&|x''| = 1\; \textrm{and}\; \delta < |x'| < \delta_c, \quad \textrm{implying}\quad
\check u - (\delta^\beta + \delta) = |x'| - (\delta^\beta + \delta) \leq |x'| - \delta,
\end{align}

\noindent
for some $c$ depending only on $\beta$ and $\delta_c$, and hence only on $n,p$.
Therefore, it follows by \eqref{eq:ner1} and \eqref{eq:svans2} that
$$
\check u - (\delta^\beta - \delta) \leq c \, u \quad \textrm{on} \quad
\partial (\{x : |x''| \leq 1\} \cap \Lambda^o_{\delta_c} \setminus \Lambda_{\delta}),
$$
for some $c(n,p)$, independent of $p$ if $p$ is large.  
By the comparison principle we obtain
$$
\check u - (\delta^\beta - \delta) \leq c \, u \quad \textrm{in} \quad
\{x : |x''| \leq 1\} \cap \Lambda^o_{\delta_c} \setminus \Lambda_{\delta},
$$
and hence, by the definition of $\check u$, 
\begin{align}\label{eq:sistaner}
|x'|^{\beta} - \delta^{\beta}  \leq  (1-|x''|^2)|x'|^{\beta} + |x'| - (\delta^\beta + \delta)  \leq  c \, u
\end{align}

\noindent
whenever $x \in \{x : |x''| = 0\} \cap \Lambda^o_{\delta_c} \setminus \Lambda_{\delta}$.

By translating the function $\check u - (\delta^\beta + \delta)$ and the domain
$\{x : |x''| \leq 1\} \cap \Lambda^o_{\delta_c} \setminus \Lambda_{\delta}$
as in the proof of the upper bound, we obtain \eqref{eq:sistaner} in $B(0,\delta_c)$.
This completes the proof of \eqref{eq:scaled_simple_0} and hence also the proof of Theorem \ref{le:pipelemma}. $\hfill \Box$ \\


\section{Estimates of $p$-harmonic measures and theorems of \\ Phragm\'en-Lindel\"of type}
\label{sec:p-harmonic-measures}
\setcounter{equation}{0} \setcounter{theorem}{0}

We first state and prove our results concerning $p$-harmonic measures.
Using these results, we then conclude our Phragm\'en-Lindel\"of-type theorems for $p$-subharmonic and $p$-harmonic functions.

In the complex plane, the harmonic measure of the semicircle
$|z| = r$, $Im(z) \geq 0$,
taken with respect to $|z| < r$, $Im(z) > 0$, is given explicitly by
\begin{align*}
v_r(z) = 2 \left(1  - \frac{1}{\pi} \arg\frac{z-r}{z+r}\right)
     =  \frac{4}{\pi} \int_0^{[z;r]^{1/4}} \frac{t\, dt}{\sqrt{1-t^4}},
\end{align*}
where $[z;r] = 4 r^2 y^2 / (4 r^2 y^2 + (r^2 - |z|^2)^2)$ and $z = x +  {\bf i}\, y$,
see e.g., Nevalinna \cite[Page 43]{Neva36} or Lindqvist \cite[Page 310]{L85}.
In $n$-dimensional space $\mathbf{R}^n$, an explicit formula is still valid in the borderline case $p = n$.
In particular, \cite[Lemma 3.5]{L85} proves the following.
Let $m \in [1, n-1]$, $\Lambda = \{ x\in \mathbf{R}^n : |x'| = 0\}$
and denote
\begin{align*}
\kappa(n,m) = \int_0^1 t^{(2 m + 1 - n)/(n - 1)} (1 - t^4)^{-1/2} dt.
\end{align*}
Define
\begin{align*}
v_r(x) = \frac{1}{\kappa(n,m)} \int_0^{[x;r]^{1/4}} t^{(2 m + 1 - n)/(n - 1)} (1 - t^4)^{-1/2} dt,
\end{align*}
where $r > 0$ and $[x;r] = {4 r^2 |x'|^2}/({4 r^2 |x'|^2 + (r^2 - |x|^2)^2})$ for $|x''|^2 \neq r^2$.
Then, $v_r(x)$ is the $n$-harmonic measure of $\partial B(0,r) \setminus \Lambda$ at $x$ 
with respect to $B(0, r)\setminus \Lambda$.
The asymptotic behaviour
\begin{align*}
C^{-1} \leq v_r(x) \, r^{m/(n-1)} \leq C
\end{align*}
as $r \to \infty$ follows, see \cite[Lemma 3.6]{L85}.

To the authors knowledge,
no explicit formula is known in the general case $p \in (n-m, \infty]$, $p \neq n$.
Nevertheless, in the below theorem, which we state and prove in more general geometry,
we show that the asymptotic behaviour,
as $r \to \infty$,
generalizes to $p \in (n-m, \infty]$ as follows.

\begin{theorem}\label{th:p-harmonic-measure}
Suppose that $m, n \in N$ such that $m \in [0,n-1]$,
let $\Lambda \subset \mathbf{R}^n$ be an $m$-dimensional hyperplane, $w \in \Lambda$, $p \in (n-m,\infty]$
and suppose that $\beta = (p-n+m)/(p-1)$ with $\beta = 1$ if $p = \infty$.
Assume that $\Omega \subset \mathbf{R}^n$ is an unbounded domain satisfying $\Lambda \subseteq \complement\Omega \subseteq \Lambda_s$ for some $s > 0$.
Let $v_r$ be the $p$-harmonic measure of $\partial B(w,5r) \setminus \complement\Omega$
with respect to $B(w, 5r)\cap \Omega$.
Then there exists $c(n,p)$, independent of $p$ if $p$ is large, such that
\begin{align*}
c^{-1}\, s^\beta \leq v_r(A_{2 s}(w))\, r^\beta \leq c \,s^\beta  
\end{align*}
whenever $2s/\delta_c < r$, where $\delta_c$ is from Lemma \ref{le:superlosningen}.
\end{theorem}

Before we prove the theorem, we make the following remark, which proof is immediate.

\begin{remark}\label{re:p-h-meas}
Using Harnack's inequality, 
Theorem \ref{th:p-harmonic-measure} implies that for any $x \in \Omega$
there exists a constant $C$ such that
\begin{align*}
C^{-1} \leq v_r(x)\, r^\beta \leq C  
\end{align*}
whenever $r$ is so large that $x \in B(w,5r)$ and $2 s/\delta_c < r$.
Moreover,
the lower bound in Theorem \ref{th:p-harmonic-measure} holds for any domain $\Omega \subset \mathbf{R}^n$
such that $\complement \Omega \subseteq \Lambda_s$. 
\end{remark}

\noindent
{\bf Proof of Theorem \ref{th:p-harmonic-measure}.}
In the following,
if $m = n-1$ so that $\Lambda$ splits $\mathbf{R}^n$ in two halves,
we focus on the upper of these halves.
To prove the upper bound,
let $\hat v$ be the $p$-harmonic function in $B(w,5r)\setminus \Lambda$,
satisfying boundary values 1 on $\partial B(w,5r)$ and $0$ on $B(w,4r) \cap \Lambda$ continuously.
If $m \geq 1$ then we also let $\hat v$ increase continuously from 0 to 1 on the set
$\Lambda \cap (B(w,5r) \setminus B(w,4r))$.
Note that the boundary values for $\hat v$ are continuous and that existence of
$\hat v$ follows from \eqref{eq:cap-estimate} and standard existence theorems,
see e.g. Heinonen--Kilpel\"ainen--Martio \cite{HKM}.
By construction of $\hat v$ and by the definition of $p$-harmonic measure (Definition \ref{def:p-hmeas}) we obtain
$$
v_r \leq \hat v \quad \textrm{in} \quad B(w,5r)\cap \Omega.
$$
%
%
%

Since $v_r(A_{5r}(w)) = 1$, we obtain by a well known H\"older continuity of the $p$-harmonic function $v_r$,
up to $\partial B(w,5r)$ near $A_{5r}(w)$ (see e.g. \cite[Lemma 2.3]{lower_order} and \cite[Lemma 2.4]{LuN}), that
$
\frac12 \leq v_r(A_{(5-\epsilon)r}(w)) \leq 1
$
for some small $\epsilon > 0$ depending only on $n$ and $p$.
Harnack's inequality now yields $c(n,p)$ such that
\begin{align}\label{eq:extra_for_reviewer}
c^{-1} \leq v_r(A_{r}(w)) \leq 1.
\end{align}
The derivation of \eqref{eq:extra_for_reviewer} is valid for the $p$-harmonic function $\hat v$ as well.
Therefore, we conclude that
$$
v_r(A_r(w)) \approx \hat v(A_r(w)) \approx 1
$$
for constants depending only on $n$ and $p$.

We next apply Theorem \ref{le:pipelemma} to $\hat v$, with $x = A_{2 s}(w)$ and $\delta = 0$,
giving
\begin{align*}
c^{-1} \, v_r(A_{2 s}(w))
\leq
\frac{\hat v(A_{2 s}(w))}{\hat v(A_r(w))}
\leq c \left(\frac{d(A_{2 s}(w), \Lambda)}{r}\right)^\beta
\leq c \, \frac{s^\beta}{r^\beta}
\end{align*}
whenever $2 s < \delta_c r$ and $c=c(n,p)$, independent of $p$ when $p$ is large.
This proves the upper bound in Theorem \ref{th:p-harmonic-measure}.

To prove the lower bound,
let $\check v$  be the $p$-harmonic function in $B(w,5r)\setminus \Lambda_{s}$, satisfying boundary values
1 on $\partial B(w,5r)\setminus \Lambda_{2 s}$ and
$0$ on $\overline{B(w,5r)}\cap \partial\Lambda_{s}$ continuously.
If $m \geq 1$ then we also let $\check v$ increase continuously from 0 to 1 on the set
$\partial B(w,5r) \cap (\Lambda_{2 s} \setminus \Lambda_{s})$.
By similar reasoning as in the proof of the upper bound we have
$$
\check v \leq v_r \quad \textrm{in} \quad B(w,5r) \setminus \Lambda_{s}
$$
and
$$
v_r(A_r(w)) \approx \check v(A_r(w)) \approx 1.
$$
We now apply Theorem \ref{le:pipelemma} to $\check v$, with $x = A_{2 s}(w)$ and $\delta = s/r$, to obtain
\begin{align*}
c^{-1} \frac{s^\beta}{r^\beta}
\leq
c^{-1} \left\{ \frac{(2 s)^\beta}{r^\beta} - \frac{s^\beta}{r^\beta}\right\}
\leq
c^{-1} \left\{ \left(\frac{d(A_{2 s}(w), \Lambda)}{r}\right)^\beta - \frac{s^\beta}{r^\beta}\right\}
\leq
\frac{\check v(A_{2 s}(w))}{\check v(A_r(w))}
\leq
c \, v_r(A_{2 s}(w))
\end{align*}
whenever $2 s < \delta_c r$ and $c=c(n,p)$, independent of $p$ when $p$ is large.
This proves the lower bound of $v_r$ and hence the proof of Theorem \ref{th:p-harmonic-measure} is complete. $\hfill \Box$ \\



We continue this section by using the estimates for $p$-harmonic measure,
given in Theorem \ref{th:p-harmonic-measure},
to prove a result of Phragmen-Lindel\"of type.
Before stating the theorem,
let us recall the classical result of Phragm\'en-Lindel\"of \cite{PL08}:
If $u(z)$, $z = x + {\bf i}\,y$, is subharmonic in the upper half plane $Im (z) > 0$,
and if $\limsup u(z) \leq 0$ as $z$ approaches any point on the real axis,
then, either $u \leq 0$ in the whole upper plane or $u$ grows so fast that
\begin{align*}
\liminf_{R \to \infty} \, \frac{\sup_{|z| = R} u(z)}{R}\, > \, 0.
\end{align*}
In the below corollary, we expand this theorem to $p$-subharmonic functions, $p \in (n-m, \infty]$,
in domains in $\mathbf{R}^n$ lying outside an $m$-dimensional hyperplane.
We note that the borderline case $p = n$ was proved by Lindqvist \cite[Theorem 4.8 and Remark 4.9]{L85},
using the explicit formula for $n$-harmonic measure, stated above Theorem \ref{th:p-harmonic-measure}.

To formulate and prove our corollary we use the notation
\begin{align*}
M(R) = \sup_{\partial B(w,R) \cap \Omega} u.
\end{align*}
\begin{corollary}\label{th:Phragmen-Lindelof}
Suppose that $m, n \in N$ such that $m \in [0,n-1]$,
let $\Lambda \subset \mathbf{R}^n$ be an $m$-dimensional hyperplane, $w \in \Lambda$, $p \in (n-m,\infty]$
and suppose that $\beta = (p-n+m)/(p-1)$ with $\beta = 1$ if $p = \infty$.
Let $\Omega$ be an unbounded domain so that $\Lambda \cap \Omega = \emptyset$.
Suppose that $u$ is $p$-subharmonic in $\Omega$
and that
\begin{align*}
\limsup_{x \to y} \, u(x) \, \leq \,  0 \quad \textrm{for all} \quad y \in \partial \Omega.
\end{align*}
Then, either $u \leq 0$ in $\Omega$ or
\begin{align*}
\liminf_{R\to\infty} \, \frac{M(R)}{R^\beta} \, > \, 0.
\end{align*}
\end{corollary}

\begin{remark}
If $\Omega = \mathbf{R}^n \setminus \Lambda_{s}$ in Corollary \ref{th:Phragmen-Lindelof},
for some $s \geq 0$,
then the $p$-harmonic function
$$
d(x,\Lambda)^\beta - s^\beta
$$
shows that the growth estimate in Corollary \ref{th:Phragmen-Lindelof} is sharp.
\end{remark}

\noindent
{\bf Proof of Corollary \ref{th:Phragmen-Lindelof}.}
The following argument is standard, see e.g. Lindqvist \cite[Principle 4.3]{L85}
or Heinonen--Kilpel\"ainen--Martio \cite[Section 11.11]{HKM}.
Assume that $u(x_0) > 0$ for some $x_0 \in \Omega$.
By the maximum principle for $p$-subharmonic functions we obtain
\begin{align*}
M(R) = \sup_{\partial B(w,R) \cap \Omega} u(x) = \sup_{B(w,R) \cap \Omega} u(x).
\end{align*}
%
Let $v_R$ be the $p$-harmonic measure for $\partial B(w,R) \setminus \Lambda$,
taken with respect to $B(w, R) \setminus \Lambda$.
Existence follows from \cite[Theorem 9.2]{HKM}.
Then
\begin{align*}
\limsup_{x \to z} u(x) \leq M(R) \,v_R(z) \quad \textrm{for all} \quad z \in \partial \left(B(w,R) \cap \Omega \right),
\end{align*}
and the comparison principle implies that $u \leq M(R) v_R$ in $B(w,R) \cap \Omega$.
Using Theorem \ref{th:p-harmonic-measure} and Remark \ref{re:p-h-meas} we have,
for any $x \in \Omega$,
the existence of a constant $C$ such that
$u(x) \leq C M(R) R^{-\beta}$ whenever $R$ is so large that $x \in B(w,R)$.
Therefore
\begin{align}\label{eq:phragmen-1}
0 < u(x_0) \leq C \frac{M(R)}{R^\beta}
\end{align}
which proves the result. $\hfill \Box$\\


We finally state and prove,
using a similar approach as in the proofs of Theorem
\ref{th:p-harmonic-measure} and Corollary \ref{th:Phragmen-Lindelof},
the following growth estimates for $p$-harmonic functions in unbounded domains:

\begin{theorem}\label{th:Phragmen-Lindelof-2}
Suppose that $m, n \in N$ such that $m \in [0,n-1]$,
let $\Lambda \subset \mathbf{R}^n$ be an $m$-dimensional hyperplane, $w \in \Lambda$, $p \in (n-m,\infty]$
and suppose that $\beta = (p-n+m)/(p-1)$ with $\beta = 1$ if $p = \infty$.
Assume that $\Omega \subset \mathbf{R}^n$ is an unbounded domain satisfying $\Lambda \subseteq \complement\Omega \subseteq \Lambda_s$ for some $s > 0$.
Suppose that $u$ is a positive $p$-harmonic function in $\Omega$,
satisfying $u = 0$ continuously on $\partial \Omega$.
Then there exists $c(n,p)$, independent of $p$ when $p$ is large, such that
\begin{align*}
c^{-1}\, s^{-\beta} d(x,\Lambda)^\beta \leq \frac{u(x)}{u(A_{2 s}(w))} \leq c \,s^{-\beta} d(x,\Lambda)^\beta
\end{align*}
whenever $x \in \mathbf{R}^n \setminus \Lambda_{2 s}$.
\end{theorem}


Theorem \ref{th:Phragmen-Lindelof-2} generalizes parts of Kilpel\"ainen--Shahgholian--Zhong \cite{KSZ07}
to more general geometries.
In particular, in \cite[Lemma 3.2]{KSZ07} it is proved that if $u$ is a non-negative $p$-harmonic function
on $\mathbf{R}^n \setminus \Lambda$, where $\Lambda$ is an $(n-1)$-dimensional hyperplane,
with $u = 0$ continuously on $\Lambda$, then
%
$
u(x) = O(|x|) 
$
%
as $|x| \to \infty$.
In the special case of $\Omega = \mathbf{R}^n\setminus\Lambda$ in Theorem \ref{th:Phragmen-Lindelof-2},
where $\Lambda$ is an $m$-dimensional hyperplane, $m \in [0, n-1]$,
we obtain the result
$
u(x) \approx d(x,\Lambda)^\beta
$
whenever $x \in \mathbf{R}^n$. \\

\noindent
{\bf Proof of Theorem \ref{th:Phragmen-Lindelof-2}.}
We begin with the lower bound.
Proceeding as in the proof of Theorem \ref{th:Phragmen-Lindelof} (from beginning to \eqref{eq:phragmen-1}) we obtain,
in place of \eqref{eq:phragmen-1},
\begin{align}\label{eq:som-i-forra-satsen}
u(A_{2 s} (w)) \leq c \,  s^\beta \, \frac{\sup_{B(w,5 r) \cap \Omega} u}{r^\beta} 
\end{align}
whenever $2 s/\delta_c < r$.
Let $\tilde c$ be the constant from Carleson's estimate
and define $\tilde u$ as the $p$-harmonic function in $B(w,5\tilde c r)\setminus \Lambda$, satisfying the boundary values
$u$ on $\partial B(w,5\tilde c r)\cap \Omega$ and
$0$ on $(\partial B(w,5\tilde c r)\setminus \Omega) \cup \Lambda$ continuously.
Then $u \leq \tilde u$ in $B(0,5\tilde c r)$ by the comparison principle and $\tilde u$
satisfies the assumptions of the Carleson estimate.
Applying the Carleson estimate to $\tilde{u}$, Harnack's inequality and a well known H\"older continuity
of $u$ and $\tilde u$ up to $\partial B(w,5\tilde c r)$ near $A_{5\tilde c r}(w)$, we obtain 
\begin{align}\label{eq:max-carlesson}
\sup_{B(w,5 r) \cap \Omega} u \,
\leq \sup_{B(w,5 r) \cap \Omega} \tilde u \,
\leq \, \tilde c \, \tilde u(A_{5 r}(w)) \,
\leq \, c \, u(A_{r}(w)).
\end{align}
By \eqref{eq:som-i-forra-satsen} and \eqref{eq:max-carlesson} we conclude that
\begin{align}\label{eq:r_och_a_r_lower}
u(A_{2 s} (w)) \leq c \, s^\beta \,  \frac{u(A_r(w))}{r^\beta}
\end{align}
whenever $2 s/\delta_c < r$, for $c = c(n,p)$, independent of $p$ if $p$ is large.
Let $\check v$ be the $p$-harmonic function in $B(w,4r)\setminus \Lambda_{s}$ satisfying the boundary values
$u$ on $\partial B(w,4r)\setminus \Lambda_{2 s}$ and $0$ on $\overline{B(w,4r)}\cap \partial\Lambda_{s}$ continuously.
If $m \geq 1$ then we also let $\check v$ increase continuously from 0 to $u$ on the set
$\partial B(w,4r) \cap (\Lambda_{2 s} \setminus \Lambda_{s})$.
By the comparison principle we obtain $\check v \leq u$ in $B(w,4r) \setminus \Lambda_{s}$.
Moreover, using Harnack's inequality and H\"older continuity we obtain
$u(A_r(w)) \approx \check v(A_r(w))$ for constants depending only on $n$ and $p$.
Applying Theorem \ref{le:pipelemma} to $\check v$ with $\delta = s/r$ gives $c = c(n,p)$ so that
\begin{align}\label{eq:johej_lower}
c^{-1} \, \left\{\left(\frac{d(x, \Lambda)}{r}\right)^{\beta} - \left(\frac{s}{r}\right)^{\beta} \right\}
\leq \frac{\check v(x)}{\check v(A_{r}(w))}
\leq c \, \frac{u(x)}{u(A_{r}(w))}
\end{align}
whenever $2 s/\delta_c < r$ and $x \in B(w,\delta_c\, r) \setminus \Lambda_{s}$.
Using \eqref{eq:r_och_a_r_lower} and \eqref{eq:johej_lower} we obtain
\begin{align*}
\frac{1}{c\, s^\beta} \, d(x, \Lambda)^\beta
\leq \frac{1}{c \, s^\beta} \left\{ d(x, \Lambda)^\beta - s^\beta \right\}
\leq \frac{u(x)}{u(A_{2 s}(w))}
\end{align*}
whenever $x \in B(w,\delta_c\, r) \setminus \Lambda_{2 s}$ and $c(n,p)$, independent of $p$ when $p$ is large.
Sending $r \to \infty$ gives the lower bound in Theorem \ref{th:Phragmen-Lindelof-2}.

To prove the upper bound, 
put $x = A_{2 s}(w)$ in \eqref{eq:johej_lower} to obtain

\begin{align}\label{eq:over-lower-and-over}
c^{-1} \frac{s^\beta}{r^\beta}
\leq
c^{-1} \left\{ \frac{(2 s)^\beta}{r^\beta} - \frac{s^\beta}{r^\beta}\right\}
\leq
\frac{\check v(A_{2 s}(w))}{\check v(A_r(w))}
\leq
c \, \frac{u(A_{2 s}(w))}{u(A_r(w))}
\end{align}
whenever $2 s / \delta_c <  r$. 
Let $\hat v$ be the $p$-harmonic function in $B(w,5r)\setminus \Lambda$ satisfying boundary values
$\sup_{B(w,5r)\cap \Omega} u$ on $\partial B(w,5r)$ and
$0$ on $B(w,4r) \cap \Lambda$ continuously.
If $m \geq 1$ then we also let $\hat v$ increase continuously from 0 to $\sup_{B(w,5r)\cap \Omega} u$ on the set
$\Lambda \cap (B(w,5r) \setminus B(w,4r))$.
By the comparison principle we obtain $u \leq \hat v$ in $B(w,5r)\cap \Omega$.
Moreover, using similar reasoning as in \eqref{eq:max-carlesson} we obtain
$u(A_r(w)) \approx \hat v(A_r(w))$ for constants depending only on $n$ and $p$.
Another application of Theorem \ref{le:pipelemma}, to $\hat v$ with $\delta = 0$,
gives
\begin{align}\label{eq:johej}
c^{-1}\, \frac{u(x)}{u(A_r(w))}
\leq
\frac{\hat v(x)}{\hat v(A_r(w))}
\leq c \, \frac{d(x,\Lambda)^{\beta}}{r^\beta}
\end{align}
whenever $x \in B(w,\delta_c\, r) \cap \Omega$ and $2 s/\delta_c < r$.
Using \eqref{eq:over-lower-and-over} and \eqref{eq:johej} we obtain
\begin{align*}
\frac{u(x)}{u(A_{2 s}(w))} \leq c \, s^{-\beta} \, d(x, \Lambda)^\beta
\end{align*}
whenever $x \in B(w,\delta_c\, r) \cap \Omega$ and $c = c(n,p)$, independent of $p$ when $p$ is large.
Sending $r \to \infty$ gives the upper bound and hence the proof of Theorem \ref{th:Phragmen-Lindelof-2} is complete. $\hfill \Box$





\bibliographystyle{amsalpha}

\end{document}